\documentclass{article}
\usepackage{amssymb}

\newtheorem{theorem}{Theorem}

\newtheorem{proposition}[theorem]{Proposition}

\input{tcilatex}
\begin{document}

\begin{center}
{\Huge Typical dynamics of Newton's method}

\bigskip

{\Large T. H. Steele}

\bigskip
\end{center}

\textbf{Abstract}: Let $C^{1}(M)$ be the set of of continuously
differentiable real valued functions defined on $[-M,M]$. We show that for
the typical element $f$ in $C^{1}(M)$, there exists a set $S\subset \lbrack
-M,M]$, both residual and of full measure in $[-M,M]$, such that for any $%
x\in S$, the trajectory $\tau (x,f)$ generated by Newton's method using $f$
and $x$ either diverges, or generates a Cantor set as its attractor.
Whenever the Cantor is an attractor, the dynamics on the attractor are
described by a single type of adding machine, so that the dynamics on all of
these attractive Cantor sets are topologically conjugate.

\textbf{Mathematics Subject Classification}: 37B20, 26A18, 54H25, 65P40

\textbf{Keywords}: Newton's method, adding machine, odometer, generic

\section{Introduction}

For a polynomial of the form $f(x)=ax^{2}+bx+c$, the quadratic formula
provides a solution to the equation $f(x)=0$. Appropriate formulae also
provide solutions to third and fourth degree equations. If a polynomial $f$
is of degree five or higher, however, no such formulae exist. This is also
the case for transcendental equations. In many calculus courses, Newton's
method is introduced as an application of the derivative that allows us to
at least approximate solutions to the equation $f(x)=0$, should $f$ be
differentiable. With Newton's method, one begins with an estimate $x_{0}$ of
a desired root $r$, and the assumption that the line tangent to $f$ at $%
x_{0} $ is close to the graph of $f$ on an interval that contains $r$. The
geometry of this situation suggests that the $x$-intercept $x_{1}$ of the
tangent line provides a better approximation of $r$. More precisely, one
takes $x_{1}$ such that $0-f(x_{0})=f^{\prime }(x_{0})(x_{1}-x_{0})$, or
that $x_{1}=\mathbf{n}(f,x_{0})=x_{0}-\frac{f(x_{0})}{f^{\prime }(x_{0})}$.
This gives rise to a sequence of successive approximations $\{x_{n}\}$ of $r$%
, given by $x_{n+1}=\mathbf{n}(f,x_{n})$. Ideally, one has $%
\lim_{n\rightarrow \infty }x_{n}=r$.

Here, we study the behavior of the trajectories $\{x_{n}\}$ generated by
Newton's method when using a continuously differentiable function $f$. We
find that rarely does the sequence $\{x_{n}\}$ converge. For any $M>0$,
there exists a residual set $\mathcal{G}$ contained in the set of
continuously differentiable functions defined on $[-M,M]$ such that for each 
$f\in \mathcal{G}$, there exists a large set $S\subset \lbrack -M,M]$ for
which the sequence generated by Newton's method does not converge whenever
the initial estimate $x_{0}$ is in $S$. More precisely, the set $S$ is both
residual and of full measure in $[-M,M]$. If $x_{0}\in S$, then one of two
possibilities occurs for the sequence $\{x_{n}\}$ generated by $x_{n+1}=%
\mathbf{n}(f,x_{n})$. The first is that for some $k\in 
\mathbb{N}
$, one has $\mid x_{k}\mid >M$. The second is that the sequence $\{x_{n}\}$
gives rise to a particular, unique type of adding machine, and generates a
Cantor set as its attractor.

\section{Preliminaries}

We work in three metric spaces. Let $C([-M,M],%
\mathbb{R}
)$ be the set of real valued continuous functions defined on $[-M,M]$, and
endow $C([-M,M],%
\mathbb{R}
)$ with the supremum metric $\Vert f-g\Vert =\sup \{\mid f(x)-g(x)\mid :x\in
I\}$. Most of our analysis takes place in $C^{1}(M)$, the set of
continuously differentiable functions defined on $[-M,M]$. Here we use the
metric $d_{1}(f,g)=\Vert f-g\Vert +\Vert f^{\prime }-g^{\prime }\Vert $.
Endowed with these metrics, both $C([-M,M],%
\mathbb{R}
)$ and $C^{1}(M)$ become complete metric spaces \cite{BBT}. Within $[-M,M]$,
we use the usual Euclidean metric. Regardless of the space we are
considering, $B_{\varepsilon }(\circ )$ represents the open ball of radius $%
\varepsilon >0$ centered at $\circ $; the nature of $\circ $ will determine
which space we are considering. Let $\lambda S$ represent the Lebesgue
measure of the set $S\subset \lbrack -M,M]$, and take $<a,b>$ to be the
closed interval with endpoints $a$ and $b$. That is, $<a,b>=[a,b]$ if $a<b$,
and $<a,b>=[b,a]$ if $b<a$.

Let $\mathbf{n}(x,f)=x-\frac{f(x)}{f^{\prime }(x)}$, and for any integer $%
m\geq 1$, $\mathbf{n}^{m}$ denotes the $m^{th}$ iterate of $\mathbf{n}$. For
each $x$ in $[-M,M]$ and $f\in C^{1}(M)$, we call $\tau (x,f)=\{\mathbf{n}%
^{m}(x,f)\}_{m=0}^{\infty }$ the trajectory of $\mathbf{n}(x,f)$. If there
exists some $k\in 
\mathbb{N}
$ such that $\mid \mathbf{n}^{k}(x,f)\mid >M$, we say that $\tau (x,f)$
diverges. Otherwise, we take the set of all subsequential limits of $\tau
(x,f)$ to be the $\omega $-limit set of $\mathbf{n}$ generated by $(x,f)$,
and write $\omega (x,f)$. Equivalently, $\omega (x,f)=\cap _{m\geq 0}%
\overline{\cup _{k\geq m}\mathbf{n}^{k}(x,f)}$.

In Proposition 4, one finds a construction critical to the sequel. We begin
with a function \thinspace $f\in C^{1}(M)$, and then develop \thinspace $%
g\in C^{1}(M)$ that well approximates $f$, and is "almost" piecewise linear.
In particular, let $\{-M=z_{0},z_{1},...,z_{m}=M\}$ be a partition of $[-M,M]
$, and to each $z_{i}$, $1\leq i\leq m-1$, associate an open interval $%
U_{i}=(a_{i},b_{i})$ such that $z_{i}\in U_{i}$, $\overline{U_{i}}\cap 
\overline{U_{j}}=\varnothing $ whenever $i\neq j$, and $%
\sum_{i=1}^{m-1}(b_{i}-a_{i})$ is as small as we like. The function $g$ will
be linear on each of the $m$ nondegenerate closed intervals complementary to 
$\cup _{i=1}^{m-1}U_{i}$ in $[-M,M]$. We use a construction due to
Misiurewicz to determine $g$ on $\cup _{i=1}^{m-1}U_{i}$ \cite{M}. This
construction allows us to extend $g$ defined on $[-M,M]-\cup
_{i=1}^{m-1}U_{i}$ to each interval $U_{i}$ such that $g^{\prime }(x)\in
<g_{-}^{\prime }(a_{i}),g_{+}^{\prime }(b_{i})>$ for all $x\in U_{i}$, $%
g^{\prime \prime }(x)$ is always either positive or negative on each $U_{i}$%
, and $\parallel f-g\parallel $ is as small as we like. In essence,
Misiurewicz's construction allows us to smoothly round the corners that
would be found at the points $\,z_{i}$, should $g$ be piecewise linear on
each interval $(z_{i},z_{i+1})$. In fact, Misiurewicz's construction
provides a function $g\in C^{\infty }$.

We recall the following conventions with the Baire category theorem in mind.
Let $(X,\rho )$ be a metric space. A set is of the first category in $X$ if
it can be written as a countable union of nowhere dense sets in $X$;
otherwise, the set is of the second category. A set is residual if it is the
complement of a first category set, and an element of a residual subset of $%
X $ is called either a typical or a generic element of $X$. The Baire
category theorem is fundamental to existence statements in the main results:

\begin{theorem}
: If $(X,\rho )$ is a complete metric space and $R$ is residual in $X$, then 
$R$ is dense in $X$.
\end{theorem}

We next develop the notion of an adding machine. Most of what follows has
been borrowed from \cite{BK} and \cite{DDS}. Let $\alpha \in (%
\mathbb{N}
-\{1\})^{%
\mathbb{N}
}$, and set $\Delta _{\alpha }=\Pi _{i=1}^{\infty }%
\mathbb{Z}
_{\alpha (i)}$, where $%
\mathbb{Z}
_{k}=\{0,1,...,k-1\}$. Take the product topology on $\Delta _{\alpha }$.
Thus, as a topological space, $\Delta _{\alpha }$ is homeomorphic to the
Cantor space. Instead of the usual coordinate-wise addition, we add two
elements of $\Delta _{\alpha }$ with "carry over" to the right. More
precisely, if $(x_{1},x_{2},....)$ and $(y_{1},y_{2},....)$ are in $\Delta
_{\alpha }$, then $(x_{1},x_{2},....)+(y_{1},y_{2},....)=(z_{1},z_{2},....)$%
, where $z_{1}=x_{1}+y_{1}$ mod $(\alpha (1))\,\ $and, in general, $z_{i}$
is defined recursively as $z_{i}=x_{i}+y_{i}+\varepsilon _{i-1}$ mod $%
(\alpha (i))$, where $\varepsilon _{i-1}=0$ if $z_{i-1}=x_{i-1}+y_{i-1}+%
\varepsilon _{i-2}$ $<$ $\alpha (i-1)$, and $\varepsilon _{i-1}=1$,
otherwise. If we let $f_{\alpha }$ be the "+1" map, that is $f_{\alpha
}(x_{1},x_{2},....)=(x_{1},x_{2},....)+(1,0,....)$, then $(\Delta _{\alpha
},f_{\alpha })$ is a dynamical system known in various contexts as a
solenoid, adding machine or odometer. For convenience, we will sometimes
refer to  $f_{\alpha }$ alone as an adding machine, with the understanding
that we are using $(\Delta _{\alpha },f_{\alpha })$. Should $(\Delta
_{\alpha },f_{\alpha })$ be an adding machine with $x\in \Delta _{\alpha }$,
it follows that $\{f_{\alpha }^{m}(x)\}_{m=0}^{\infty }$ is dense in $\Delta
_{\alpha }$, and that $\Delta _{\alpha }$ is a minimal set.

Fix $\alpha \in (%
\mathbb{N}
-\{1\})^{%
\mathbb{N}
}$, and define a function $M_{\alpha }$ from the set of primes into $%
\{0\}\cup 
\mathbb{N}
\cup \{\infty \}$ so that for each prime $p$, one takes $M_{\alpha
}(p)=\Sigma _{i=1}^{\infty }n(i)$, where $n(i)$ is the largest power of $p$
which divides $\alpha (i)$. The following theorem is a beautiful
characterization of adding machines up to topological conjugacy \cite{BS}, 
\cite{BK}.

\begin{theorem}
: Let $\alpha ,\beta \in (%
\mathbb{N}
-\{1\})^{%
\mathbb{N}
}$. Then $f_{\alpha }$ and $f_{\beta }$ are topologically conjugate if and
only if $M_{\alpha }=M_{\beta }$.
\end{theorem}

The following useful theorem is from \cite{BK}.

\begin{theorem}
: Let $\alpha \in (%
\mathbb{N}
-\{1\})^{%
\mathbb{N}
}$. Let $m_{i}=\alpha (1)\alpha (2)...\alpha (i)$ for each $i$. Let $%
f:X\rightarrow X$ be a continuous map of a compact topological space $X$.
Then $f$ is topologically conjugate to $f_{\alpha }$ if and only if the
following hold:
\end{theorem}

\begin{enumerate}
\item \textit{For each positive integer }$i$\textit{, there is a cover }$\Pi
_{i}$\textit{\ of }$X$\textit{\ consisting of }$m_{i}$\textit{\ pairwise
disjoint, nonempty, clopen sets which are cyclically permuted by }$f$\textit{%
.}

\item \textit{For each positive integer }$i$\textit{, }$\Pi _{i+1}$\textit{\
partitions }$\Pi _{i}$\textit{.}

\item \textit{If }$W_{1}\supset W_{2}\supset W_{3}\supset ....$\textit{\ is
a nested sequence with }$W_{i}\in $\textit{\ }$\Pi _{i}$\textit{\ for each }$%
i$\textit{, then }$\cap _{i=1}^{\infty }W_{i}$\textit{\ consists of a single
point.}
\end{enumerate}

Of particular interest in what follows are $\infty $-adic adding machines.
These are adding machines associated with those $\alpha $ for which $%
M_{\alpha }(p)=\infty $ for all prime numbers $p$. Thus, for each prime
number $p$, there exist infinitely many indices $i$ for which $p$ divides $%
m_{i}$, the cyclic period of the covering sets $\{W_{1},W_{2},...,W_{m_{i}}\}
$.

\section{Results on $C^{1}(M)$}

\begin{proposition}
: Let $x^{\prime }\in (-M,M)$, $f\in C^{1}(M)$ and $\varepsilon >0$. There
exists $g\in C^{1}(M)$ such that $d_{1}(f,g)<\varepsilon $, and either $\tau
(x^{\prime },g)$ is divergent, or $\omega (x^{\prime },g)$ is periodic.
\end{proposition}

\textit{Proof: }We work in $C^{1}(M)$. Let $f\in C^{1}(M)$, $x^{\prime }\in
(-M,M)$ and $\varepsilon >0$. Since $f\in C^{1}(M)$, there exists $\delta
_{1}>0$ such that $\mid f(x)-f(y)\mid <\frac{\varepsilon }{3}$ whenever $%
\mid x-y\mid <$ $\delta _{1}$, and there exists $\delta _{2}>0$ such that $%
\mid f^{\prime }(x)-f^{\prime }(y)\mid <\frac{\varepsilon }{3}$ whenever $%
\mid x-y\mid <$ $\delta _{2}$. Let $\delta =\min \{\delta _{1},\delta _{2}\}$%
. Take $\{-M=z_{0},z_{1},...,z_{m}=M\}$ to be a $\delta $-fine partition of $%
[-M,M]$ such that $x^{\prime }\in (z_{j},z_{j+1})$, for some $0\leq j<m$.
Consider the corresponding values $\{f(z_{0}),f(z_{1}),...,f(z_{m})\}$. We
perturb the values $f(z_{i})$ to $g_{1}(z_{i})$ so that

\begin{enumerate}
\item $g_{1}(z_{i})\neq g_{1}(z_{j})$ whenever $i\neq j$.

\item $\mid f(z_{i})-g_{1}(z_{i})\mid <\frac{\varepsilon }{6}$ for all $%
0\leq i\leq m$.

\item $\mid \lbrack \frac{f(z_{i})-f(z_{i+1})}{z_{i}-z_{i+1}}]-[\frac{%
g_{1}(z_{i})-g_{1}(z_{i+1})}{z_{i}-z_{i+1}}]\mid <\frac{\varepsilon }{6}$
for all $0\leq i<m$.

\item Say that the line with slope $\frac{g_{1}(z_{i})-g_{1}(z_{i+1})}{%
z_{i}-z_{i+1}}$ passing through the point $(z_{i},g_{1}(z_{i}))$ has $y_{i}$
as its $x$-intercept. Then $y_{i}\notin \{z_{0},z_{1},...,z_{m}\}$.
\end{enumerate}

Now, extend $g_{1}$ linearly to all of $[-M,M]$. Since $x^{\prime }\in
(z_{j},z_{j+1})$ for some $0\leq j<m$, there exists $\eta >0$ such that $%
B_{\eta }(x^{\prime })\subset (z_{j},z_{j+1})$. Now, fix $\alpha >0$. There
exist open intervals $U_{i}$, $0<i<m$, such that $\overline{U_{i}}\cap 
\overline{U_{j}}=\varnothing $ whenever $i\neq j$, $z_{i}\in U_{i}$ for all $%
0<i<m$, $\{x^{\prime },\cup _{i=0}^{m-1}y_{i}\}\cap (\cup _{i=1}^{m-1}%
\overline{U_{i}})=\varnothing $, and $\sum_{i=1}^{m-1}\mid U_{i}\mid <\alpha 
$. Let $g=g_{1}$ on $[-M,M]-\cup _{i=1}^{m-1}U_{i}$, and extend $g$ to $\cup
_{i=1}^{m-1}U_{i}$ using Misiurewicz's construction so that $%
d_{1}(f,g)<\varepsilon $.

Since $g$ is linear on each of the $m$ closed intervals $J_{j}\subset
(z_{j},z_{j+1})$ comprising $[-M,M]-\cup _{i=1}^{m-1}U_{i}$, it follows that 
$\mathbf{n}(z,g)=\mathbf{n}(y,g)$ whenever $z$ and $y$ are both contained in
some $int(J_{j})$. Moreover, $\mathbf{n}(x^{\prime },g)\subset 
\mathbb{R}
-\cup _{i=1}^{m-1}U_{i}$ for any $x^{\prime }\in $ $(-M,M)-\cup
_{i=1}^{m-1}U_{i}$. It follows that $\tau (x^{\prime },g)$ must either
diverge or be eventually periodic. $\square $

The following observations follow from Proposition 4 and the construction
found in its proof. Since they are critical to what follows, we highlight
them at this time.

\textbf{Observation 1}: Let $f\in C^{1}(M)$, $\varepsilon >0$ and $\delta >0$%
. There exist $g\in C^{1}(M)$ such that $d_{1}(f,g)<\varepsilon $ and a
closed set $S\subset \lbrack -M,M]$ composed of finitely many nondegenerate
closed intervals such that $\lambda S>2M-\delta $, disjoint open sets $D$
and $P$ such that $D\cup P=int(S)$, and

\begin{enumerate}
\item if $x\in D$, then $\tau (x,g)$ is divergent, and

\item if $x\in P$, then $\omega (x,g)$ is periodic.
\end{enumerate}

\textbf{Observation 2}: Take $g\in C^{1}(M)$ as found in Proposition 4.
Suppose that $\omega (x,g)$ is periodic, and $x=x_{0},x_{1},...x_{m}$ are
the distinct points found in $\tau (x,g)\subset int(S)$ such that $\mathbf{n}%
(x_{i},g)=x_{i+1}$ for all $0\leq i\leq m-1$, and $\mathbf{n}(x_{m},g)=x_{l}$
for some $0\leq l\leq m$.

\begin{enumerate}
\item There exist closed intervals $J_{i}$ for $0\leq i\leq m$ such that $%
x_{i}\in int(J_{i})\subset J_{i}\subset int(S)$ for each $i$, $\mathbf{n}%
(J_{i},g)=x_{i+1}$ for all $0\leq i\leq m-1$, and $\mathbf{n}(J_{m},g)=x_{l}$%
.

\item Consequently, there exists $\eta >0$ such that if $d_{1}(h,g)<\eta $,
then $\mathbf{n}(J_{i},h)\subset int(J_{i+1})$ for all $0\leq i\leq m-1$,
and $\mathbf{n}(J_{m},h)\subset int(J_{l})$.
\end{enumerate}

\textbf{Observation 3}: Take $g\in C^{1}(M)$ as found in Proposition 4.
Suppose that $\tau (x,g)$ is divergent with $x_{i}=\mathbf{n}^{i}(x,g)\in
int(S)$ for all $i<K$, and $\mid x_{K}\mid >M$.

\begin{enumerate}
\item There exist closed intervals $J_{i}$ for $0\leq i<K$ such that $%
x_{i}\in int(J_{i})\subset J_{i}\subset int(S)$ for each $i$, $\mathbf{n}%
(J_{i},g)=x_{i+1}$ for all $0\leq i\leq K-2$, and $\mid \mathbf{n}%
(J_{K-1},g)\mid >M$.

\item Consequently, there exists $\eta >0$ such that if $d_{1}(h,g)<\eta $,
then $\mathbf{n}(J_{i},h)\subset int(J_{i+1})$ for all $0\leq i\leq K-2$,
and $\mathbf{n}(J_{K-1},h)\cap \lbrack -M,M]=\varnothing $.
\end{enumerate}

With Propositions 5 and 6, we develop those periodic sets $W_{i}$ found in
Theorem 3 that give rise to the $\infty $-adic adding machines. Proposition
5 will insure that for any prime number $p$, $p$ divides $m_{i}$ infinitely
often. Proposition 6 insures that the partitions $\Pi _{i}$ generated by
Proposition 5 give rise to $\infty $-adic adding machines. 

\begin{proposition}
: Take $g\in C^{1}(M)$ as found in Proposition 4, $\varepsilon >0$ and $t$ a
natural number. Suppose that $\omega (x,g)\subset int(S)$ is periodic. Then
there exists $h\in C^{1}(M)$ as found in Proposition 4 such that $%
d_{1}(h,g)<\varepsilon $, $\omega (x,h)$ is periodic, and $t$ divides $\mid
\omega (x,h)\mid $, the period of $\omega (x,h)$.
\end{proposition}

\textit{Proof: }Let us continue with the notation established for $g$ in
Observation 2. Let $\varepsilon >0$, and suppose that $l<m$, and $l\,<i\leq
m $. The point $x_{i-1}$ is contained in an interval $J_{i-1}=[a,b]$ on
which $g$ is linear. Without loss of generality, suppose that $\mid
x_{i}-a\mid <\mid x_{i}-b\mid $, so that $J_{i-1}$ lies to the right of $%
x_{i}$. Consider the line passing through the point $(a,g(a))$ with slope $%
g^{\prime }(x_{i-1})+\frac{\varepsilon }{3}$. Say that its $x$-intercept is $%
\alpha $. Consider the line passing the point $(a,g(a))$ with slope $%
g^{\prime }(x_{i-1})-\frac{\varepsilon }{3}$. Say that its $x$-intercept is $%
\beta $. Then $x_{i}\in int(<\alpha ,\beta >)$, where $<\alpha ,\beta >$ is
the closed interval with endpoints $\alpha $ and $\beta $. Let $y\in
int([c,d])\cap int(<\alpha ,\beta >)$, where $x_{i}$ is contained in the
interval $J_{i}=[c,d]$ on which $g$ is linear. We show that if $x^{\prime }$
is any point of $int(J_{i-1})$ and $\sigma >0$, then there exists $h\in
C^{1}(M)$ as found in Proposition 4 such that $d_{1}(h,g)<\varepsilon $, $%
h=g $ on $[-M,M]-B_{\delta }(x^{\prime })$ for some $0<\delta <\sigma $, and 
$\mathbf{n}(x^{\prime },h)=y$.

Take $0<\delta <\sigma $ sufficiently small and $h$ linear on $B_{\frac{%
\delta }{3}}(x^{\prime })$ such that

\begin{enumerate}
\item $B_{\delta }(x^{\prime })\subset int(J_{i-1})$,

\item $h(x^{\prime })=g(x^{\prime })$,

\item $\mathbf{n}(x^{\prime },h)=y$, and

\item $\parallel h-g\parallel <\frac{\varepsilon }{3}$ on $B_{\frac{\delta }{%
3}}(x^{\prime })$.
\end{enumerate}

Since $y\in int(J_{i})\cap int(<\alpha ,\beta >)$, it follows that $%
\parallel h^{\prime }-g^{\prime }\parallel =\mid h^{\prime }(x^{\prime
})-g^{\prime }(x^{\prime })\mid <\frac{\varepsilon }{3}$ on $B_{\frac{\delta 
}{3}}(x^{\prime })$, and $d_{1}(h,g)<\frac{2\varepsilon }{3}$ on $B_{\frac{%
\delta }{3}}(x^{\prime })$. Now, extend $h$ to all of $B_{\delta }(x^{\prime
})$ using Misiurewicz's construction such that $d_{1}(h,g)<\varepsilon $ on
all of $B_{\delta }(x^{\prime })$ and $h=g$ on $[-M,M]-B_{\delta }(x^{\prime
})$.

Recall that $\tau (x,g)$ is eventually periodic with $%
x=x_{0},x_{1},...,x_{l},...,x_{m}$ being the distinct points of $\tau (x,g)$
such that $\mathbf{n}(x_{i},g)=x_{i+1}$ for all $0\leq i\leq m-1$, and $%
\mathbf{n}(x_{m},g)=x_{l}$. As discussed in the previous paragraph, to each $%
x_{i}$, $l\leq i\leq m$, we associate $\delta _{i}>0$. In each open ball $%
B_{\delta _{i}}(x_{i})$, we choose $t$ distinct points $x_{i}^{j}$, $1\leq
j\leq t$, such that $x_{i}=$ $x_{i}^{q}$ for some $1\leq q\leq t$. We now
take $h\in C^{1}(M)$ such that $d_{1}(h,g)<\varepsilon $, $h=g$ on $%
[-M,M]-\cup _{i=l}^{m}B_{\delta _{i}}(x_{i})$, and if $l\leq i<m$, then for
the map $\mathbf{n}(\circ ,h)$ we have $x_{i}^{j}\mapsto x_{i+1}^{j}$ for
all $1\leq j\leq t$, $x_{m}^{j}\mapsto x_{l}^{j+1}$ for all $1\leq j<t$, and 
$x_{m}^{t}\mapsto x_{l}^{1}$. It follows that $\tau (x,h)$ is eventually
periodic, as the trajectory terminates in the $[(m-l)+1]t$ cycle $%
x_{i}^{j}\mapsto x_{i+1}^{j}$ for $l\leq i<m$ and $1\leq j\leq t$, $%
x_{m}^{j}\mapsto x_{l}^{j+1}$ for $1\leq j<t$, and $x_{m}^{t}\mapsto
x_{l}^{1}$.

Now, suppose that $l=m$. Set $z=x_{l}=x_{m}$, and as before, take $[a,b]$
such that $z\in (a,b)$, and $g$ is linear on $[a,b]$. Without loss of
generality, suppose that $g^{\prime }(z)>0$, and take $\varepsilon >0$. Let $%
\varepsilon >\sigma >0$ such that $B_{\sigma }(z)\subset (a,b)$. We begin to
construct $h\in C^{1}(M)$. Take $h(z)=-\frac{\sigma }{p}$, $h^{\prime
}(z)=g^{\prime }(z)-\frac{\sigma }{p}$ and $p\geq 3$ minimal such that $%
z<x_{1}=\mathbf{n}(z,h)\,<z+\sigma $. For $i=1,2,...,t-1$, we take $%
h(x_{i-1})=g(x_{i-1})$, $h^{\prime }(x_{i-1})=g^{\prime }(x_{i-1})+\frac{%
\sigma }{2p}$, and $x_{i}=\mathbf{n}(x_{i-1},h)$. Set $x_{t}=z$, so that $%
h(x_{t-1})=g(x_{t-1})$ and $h^{\prime }(x_{t-1})=g^{\prime }(x_{t-1})$. Let $%
\delta ^{\prime }=\min \{\mid x_{i}-x_{j}\mid :i\neq j\}$, and take $%
0<\delta <\delta ^{\prime }$ such that if $h$ is linear on each $B_{\frac{%
\delta }{5}}(x_{i})$, $1\leq i\leq t-2$, with $h(x_{i})=g(x_{i})$, and $%
h^{\prime }(x_{i})=g^{\prime }(x_{i})+\frac{\sigma }{2p}$, then $\parallel
h-g\parallel <\frac{\sigma }{2p}$ on $\cup _{i=1}^{t-2}B_{\frac{\delta }{5}%
}(x_{i})$. On $B_{\frac{\delta }{5}}(x_{t-1})$, take $h$ linear with $%
h(x_{t-1})=g(x_{t-1})$, and $h^{\prime }(x_{t-1})=g^{\prime }(x_{t-1})$. On $%
B_{\frac{\delta }{5}}(z)$, take $h$ linear with $h(z)=-\frac{\sigma }{p}$,
and $h^{\prime }(z)=g^{\prime }(z)-\frac{\sigma }{p}$. Thus, $d_{1}(h,g)\,<%
\frac{2\sigma }{p}\leq \frac{2\sigma }{3}<\frac{2\varepsilon }{3}$ on $\cup
_{i=1}^{t}B_{\frac{\delta }{5}}(x_{i})$. Let $h=g$ on $[-M,M]-B_{\sigma }(z)$%
, and extend $h$ to $B_{\sigma }(z)-\cup _{i=1}^{t}B_{\frac{\delta }{5}%
}(x_{i})$ such that $d_{1}(h,g)<\frac{2\sigma }{p}$ there, too. Thus, we
have $h\in C^{1}(M)$ such that $d_{1}(h,g)\,<\varepsilon $, and $\tau (x,h)$
is eventually periodic with period $t$. $\square $

\textbf{Observation 4}: Let $\{f_{i}\}_{i=1}^{\infty }$ be dense in $%
C^{1}(M) $. For any $f_{i}$, $\varepsilon >0$ and $j\in 
\mathbb{N}
$, there exists $g_{i,j}\in C^{1}(M)$ as found in Proposition 4, and $%
S_{i,j}\subset \lbrack -M,M]$ as described in Observation 1, such that $%
d_{1}(f_{i},g_{i,j})<\frac{\varepsilon }{2^{i+j}}$, $\lambda S_{i,j}>2M-%
\frac{\varepsilon }{2^{i+j}}$ and disjoint open sets $D_{i,j}$ and $P_{i,j}$
such that $D_{i,j}\cup P_{i,j}=int(S_{i,j})$, and

\begin{enumerate}
\item if $x\in D_{i,j}$, then $\tau (x,g_{i,j})$ diverges, and

\item if $x\in P_{i,j}$, then $\omega (x,g_{i,j})$ is periodic. Moreover, $%
j! $ divides $\mid \omega (x,g_{i,j})\mid $, the period of $\omega
(x,g_{i,j})$.
\end{enumerate}

\begin{proposition}
: Let $n$ be a natural number. There exists $\mathcal{G}_{n}$ a dense $%
G_{\delta }$ subset of $C^{1}(M)$ such that for any $h\in \mathcal{G}_{n}$,
there exists $S\subset \lbrack -M,M]$ such that $\lambda S>2M-\frac{1}{n}$,
and if $x\in S$, then either $\tau (x,h)$ diverges, or $\omega (x,h)$ is an $%
\infty $-adic odometer.
\end{proposition}

\textit{Proof:} Let $\varepsilon >0$. Take $g_{i,j}\in C^{1}(M)$ as found in
Observation 4, and consider $[a,b]$ a component of the set $S_{i,j}$. If $%
\tau (x,g_{i,j})$ diverges for any $x\in \lbrack a,b]$, take $\eta _{\lbrack
a,b]}>0$ such that $d_{1}(h,g_{i,j})<\eta _{\lbrack a,b]}$ implies that $%
\tau (x,h)$ diverges, too. That this is possible follows from Observation 3.
If $\omega (x,g_{i,j})$ is periodic for any $x\in \lbrack a,b]$, take $\eta
_{\lbrack a,b]}>0$ and intervals $J_{p}$ in accordance with Observation 2
such that $\mid J_{p}\mid <\frac{\varepsilon }{2^{i+j}}$ for all $p$, and if 
$d_{1}(h,g_{i,j})<\eta _{\lbrack a,b]}$, then $\mathbf{n}^{l}([a,b],h)%
\subset int(J_{l})$, $\mathbf{n}(J_{p},h)\subset int(J_{p+1})$, for $l\leq
p\leq m-1$, and $\mathbf{n}(J_{m},h)\subset int(J_{l})$. Since $S_{i,j}$ is
composed of finitely many pairwise disjoint non-degenerate closed intervals,
and $g_{i,j}$ is linear on each of these components, the function $\mathbf{n}%
(\circ ,g_{i,j})$ generates on $S_{i,j}$ only finitely many trajectories
with distinct tails. Let $\eta =\min \{\eta _{\lbrack a,b]}:[a,b]$ is a
component of $S_{i,j}\}$, and take $0\,<\eta _{i,j}<\min \{\eta ,\frac{%
\varepsilon }{2^{i+j}}\}$. Set $\mathcal{G}_{j}=\cup _{i=1}^{\infty }B_{\eta
_{i,j}}(g_{i,j})$. Then $\mathcal{G}_{j}$ is a dense and open subset of $%
C^{1}(M)$. Set $\mathcal{G}=\cap _{i=1}^{\infty }\mathcal{G}_{j}$. Then $%
\mathcal{G}$ is a dense $G_{\delta }$ subset of $C^{1}(M)$. Take $h\in 
\mathcal{G}$. Since $h\in \mathcal{G}_{j}$, there exists $j(i)$ such that $%
h\in B_{\eta _{j(i),j}}(g_{j(i),j})$ and $S_{j(i),j}$ such that $\lambda
S_{j(i),j}>2M-\frac{\varepsilon }{2^{j(i)+j}}>2M-\frac{\varepsilon }{2^{j}}$%
. Let $S=[-M,M]-\cup _{j=1}^{\infty }S_{j(i),j}$. Then $\lambda
S>[-M,M]-\varepsilon $. If $x\in S$, then either

\begin{enumerate}
\item $x\in $ $D_{j(i),j}$ for some $j(i)$ and $j$, and by Observations 4
and 3, $\tau (x,h)$ diverges. Otherwise,

\item $x\in $ $P_{j(i),j}$ for each $j$. By Observations 2 and 4, there
exist pairwise disjoint nondegenerate closed intervals $%
J_{1}^{j},J_{2}^{j},...,J_{m}^{j}$ such that $\mathbf{n}(J_{p}^{j},h)\subset
int(J_{p+1}^{j})$ for all $0\leq p\leq m-1$, and $\mathbf{n}%
(J_{m}^{j},h)\subset int(J_{l}^{j})$, with $\omega (x,h)\subset \cup
_{p=l}^{m}J_{p}^{j}$ and $\omega (x,h)\cap J_{p}^{j}\neq \varnothing $ for
each $l\leq p\leq m$. It follows that $P_{j}=%
\{J_{l}^{j},J_{l+1}^{j},...,J_{m}^{j}\}$ partitions $\omega (x,h)$ into $%
(m-l)+1$ periodic portions, and by construction, $j!$ divides $(m-l)+1$. By
taking an appropriate subsequence $\{j_{k}\}$, we have that $P_{j_{k}+1}$
refines $P_{j_{k}}$. We conclude that $\omega (x,h)$ is an $\infty $-adic
odometer. $\square $
\end{enumerate}

With Proposition 6, we are in a position to prove our main results. Theorem
7 deals with measure, and Theorem 8 addresses category. Theorem 9 follows
readily from Theorem 8 and the Kuratowski-Ulam theorem.

\begin{theorem}
: There exists $\mathcal{G}$ a dense $G_{\delta }$ subset of $C^{1}(M)$ such
that for any $h\in \mathcal{G}$, there exists $S\subset \lbrack -M,M]$ such
that $\lambda S=2M$, and if $x\in S$, then either $\tau (x,h)$ diverges, or $%
\omega (x,h)$ is an $\infty $-adic odometer.
\end{theorem}

\textit{Proof:} Continuing with the notation established in the previous
theorem, set $\mathcal{G=\cap }_{n=1}^{\infty }\mathcal{G}_{n}$. Since each $%
\mathcal{G}_{n}$ is a dense $G_{\delta }$ subset of $[-M,M]$, so is $%
\mathcal{G}$. If $h\in \mathcal{G}$, then for any $n\in 
\mathbb{N}
$ there exists $S_{n}$ such that $\lambda S_{n}>2M-\frac{1}{n}$, and if $%
x\in S_{n}$, then either $\tau (x,h)$ diverges, or $\omega (x,h)$ is an $%
\infty $-adic odometer. Set $S=\cup _{n=1}^{\infty }S_{n}$. Then $\lambda
S=2M$, and if $x\in S$, then either $\tau (x,h)$ diverges, or $\omega (x,h)$
is an $\infty $-adic odometer. $\square $

\begin{theorem}
: There exists $\mathcal{G}$ a dense $G_{\delta }$ subset of $[-M,M]\times
C^{1}(M)$ such that for any $(x,h)\in \mathcal{G}$, either $\tau (x,h)$
diverges, or $\omega (x,h)$ is an $\infty $-adic odometer.
\end{theorem}

\textit{Proof: Let }$\{(x_{i},f_{i})\}_{i=1}^{\infty }$ be a dense subset of 
$(-M,M)\times C^{1}(M)$, and $\varepsilon >0$. Similarly to what is found in
the proof of Proposition 6, we use Observation 4 to associate to each $%
(i,j)\in 
\mathbb{N}
\times 
\mathbb{N}
$ a function $g_{i,j}$, a set $S_{i,j}\subset \lbrack -M,M]$ and $\eta
_{i,j}>0$ such that $B_{\eta _{i,j}}(x_{i})\subset int(S_{i,j})$. From
Observations 2 and 3, either

\begin{enumerate}
\item $\tau (y,h)$ diverges for all $(y,h)\in B_{\eta _{i,j}}(x_{i})\times
B_{\eta _{i,j}}(g_{i,j})$, or

\item for each $(y,h)\in B_{\eta _{i,j}}(x_{i})\times B_{\eta
_{i,j}}(g_{i,j})$, there exist closed intervals $J_{k}^{j}$, $1\leq k\leq m$
such that

\begin{enumerate}
\item $\mid J_{k}^{j}\mid \,<\frac{\varepsilon }{2^{i+j}}$, for all $k$,

\item $y\in int(J_{1}^{j})$,

\item $\mathbf{n}(J_{k}^{j},h)\subset int(J_{k+1}^{j})$, for $1\leq k\leq
m-1 $,

\item $\mathbf{n}(J_{m}^{j},h)\subset int(J_{l}^{j})$, for some $1\leq l\leq
m$, and

\item $j!$ divides $(m-l)+1$.
\end{enumerate}
\end{enumerate}

Let $\mathcal{G}_{j}=\cup _{i=1}^{\infty }(B_{\eta _{i,j}}(x_{i})\times
B_{\eta _{i,j}}(g_{i,j}))$, a dense and open subset of $(-M,M)\times C^{1}(M)
$, and take $\mathcal{G=\cap }_{j=1}^{\infty }\mathcal{G}_{j}$, so that $%
\mathcal{G}$ is a dense $G_{\delta }$ subset of $[-M,M]\times C^{1}(M)$.
Should $(y,h)\in \mathcal{G}$, and $(y,h)\in $ $D_{i,j}\times $ $B_{\eta
_{i,j}}(g_{i,j})$ for some $i$ and $j$, then $\tau (y,h)$ diverges. Fix $%
j\in 
\mathbb{N}
$, and suppose that $\tau (y,h)$ does not diverge. As in the proof of
Proposition 6, there exist intervals $P_{j}=%
\{J_{l}^{j},J_{l+1}^{j},...,J_{m}^{j}\}$ that partition $\omega (y,h)$ into $%
(m-l)+1$ periodic portions. By taking an appropriate subsequence $\{j_{k}\}$%
, we have that $P_{j_{k}+1}$ refines $P_{j_{k}}$. The considerations found
in (2) above allow us to conclude that $\omega (y,h)$ is an $\infty $-adic
odometer. $\square $

\begin{theorem}
: There exists $\mathcal{G}$ a dense $G_{\delta }$ subset of $C^{1}(M)$ such
that for any $h\in \mathcal{G}$, there exists $S\subset \lbrack -M,M]$
residual in $[-M,M]$ such that $\lambda S=2M$, and if $x\in S$, then either $%
\tau (x,h)$ diverges, or $\omega (x,h)$ is an $\infty $-adic odometer.
\end{theorem}

\textit{Proof:} Let $\mathcal{G}_{1}$ be the dense $G_{\delta }$ subset of $%
[-M,M]\times C^{1}(M)$ found in the conclusion of Theorem 8. From the
Kuratowski-Ulam Theorem \cite{O}, it follows that there exists $\mathcal{G}%
_{2}$ a dense $G_{\delta }$ subset of $C^{1}(M)$ such that for any $f\in 
\mathcal{G}_{2}$, there exists $S_{f}$ a dense $G_{\delta }$ subset of $%
[-M,M]$ such that $S_{f}\times \{f\}\subset \mathcal{G}_{1}$. Now, let $%
\mathcal{G}_{3}$ be the dense $G_{\delta }$ subset of $C^{1}(M)$ found in
the conclusion of Theorem 7, and set $\mathcal{G}=\mathcal{G}_{2}\cap 
\mathcal{G}_{3}$. If $f\in \mathcal{G}$, then there exists $S_{f}$ a dense $%
G_{\delta }$ subset of $[-M,M]$ such that either $\tau (x,f)$ diverges, or $%
\omega (x,f)$ is an $\infty $-adic odometer, whenever $x\in S_{f}$. From
Theorem 7, there exists $S_{f}^{^{\prime }}\subset \lbrack -M,M]$ such that $%
\lambda S_{f}^{^{\prime }}=2M$, and if $x\in S_{f}^{^{\prime }}$, then
either $\tau (x,f)$ diverges, or $\omega (x,f)$ is an $\infty $-adic
odometer. Set $S=S_{f}\cup S_{f}^{^{\prime }}$. It follows that $S$ contains
a dense $G_{\delta }$ subset of $[-M,M]$, $\lambda S=2M$, and if $x\in S$,
then either $\tau (x,h)$ diverges, or $\omega (x,h)$ is an $\infty $-adic
odometer. $\square $

\bigskip

T. H. Steele

Department of Mathematics

Weber State University

Ogden UT 84408-2517

USA

thsteele@weber.edu

\end{document}


\maketitle
\renewcommand{\thefootnote}{}
\footnote{\emph{Mathematics Subject Classification}: 37B20, 26A18, 54H25, 65P40}
\footnote{\emph{Key words}: Newton's method, Adding machine, Typical behavior}
\renewcommand{\thefootnote}{\arabic{footnote}}
\setcounter{footnote}{0}

\begin{abstract}
Let $C^1(M)$ be the space of continuously differentiable real-valued functions defined on $[-M,M]$. Here, we address an irremediable flaw found in \cite{Steele}, and show that for the typical element $f$ in $C^1(M)$, there exists a set $S \subseteq [-M,M]$, both residual and of full measure in $[-M,M]$, such that for any $x \in S$, the trajectory generated by Newton's method using $f$ and $x$ either diverges, converges to a root of $f$, or generates a Cantor set as its attractor. Whenever the Cantor set is the attractor, the dynamics on the attractor are described by a single type of adding machine, so that the dynamics on all of these attracting Cantor sets are topologically equivalent.
\end{abstract}

\section{Introduction}
In the course of preparing a seminar for the Faculty of Mathematics and Physics at Charles University, Jan Dud\'ak found an irremediable flaw in the proof of \cite[Proposition 5]{Steele}. The purpose of this submission is not only to announce the error, but also to present a complete and correct treatment of the dynamics of Newton's method for continuously differentiable functions.

The error found in the proof of \cite[Proposition 5]{Steele} comes at its conclusion. There, Steele asserts that whenever $z$ is a root of some $f \in C^1(M)$ that is affine on an interval $(a,b)$ containing $z$, then there exists $g \in C^1(M)$, arbitrarily close to $f$, and $x \in (a,b)$, such that the trajectory generated by Newton's method using $g$ and $x$ is eventually periodic, and that the period is greater than one. This assertion is clearly negated by Proposition \ref{PropositionConvergenceOfNewtonsMethod}, that was developed by Dud\'ak.

\section{Preliminaries and notation}

We denote the set of all natural numbers (i.e. positive integers) by $\N$. Moreover, we denote by $\N_0$ the set of all nonnegative integers. In other words, $\N_0 = \N \cup \{ 0 \}$.

For any subset $A$ of a topological space $X$, we denote by $\overline{A}$ and $\Int (A)$ the closure of $A$ and the interior of $A$, respectively. A subset of a topological space $X$ is said to be \textit{clopen} if it is both open and closed in $X$. A partition $\mathcal{P}$ of a topological space $X$ is said to be \textit{clopen} if each set $A \in \mathcal{P}$ is clopen in $X$.

A subset $A$ of a topological space $X$ is said to be \textit{meagre} in $X$ if it is the union of countably many nowhere dense subsets of $X$. Moreover, $A$ is \textit{comeagre} in $X$ if $X \setminus A$ is meagre in $X$.

Recall that the Baire category theorem states that the intersection of countably many dense open subsets of a complete metric space is dense. Consequently, a subset $A$ of a complete metric space $X$ is comeagre in $X$ if and only if there is a dense $G_\delta$ set $D$ in $X$ such that $D \subseteq A$.

A \textit{topological dynamical system} is a pair $(X,f)$, where $X$ is a topological space and $f \colon X \to X$ is a continuous mapping. Topological dynamical systems $(X,f)$, $(Y,g)$ are said to be \textit{topologically conjugate} if there is a homeomorphism $h \colon X \to Y$ such that $g = h \circ f \circ h^{-1}$.

A sequence $(x_n)_{n=0}^\infty$ is said to be \textit{eventually periodic} if there is $p \in \N$ such that $x_{n+p}=x_n$ for all but finitely many $n \in \N_0$. If this is the case, the least such $p$ is called the period of the sequence. If the period is equal to $1$, we say that the sequence is \textit{eventually constant}.

We denote by $\lambda$ the Lebesgue measure on $\R$.

\begin{definition}
    For every $\alpha \in (\N \setminus \{ 1 \})^\N$, we denote $\Delta_\alpha := \prod_{i=1}^\infty \mathbb{Z}_{\alpha (i)}$ and equip $\Delta_\alpha$ with the product topology. Define a mapping $\varphi_\alpha \colon \Delta_\alpha \to \Delta_\alpha$ as follows. For any $x=(x_i)_{i=1}^\infty \in \Delta_\alpha$, let $\varphi_\alpha (x)$ be the member $(z_i)_{i=1}^\infty$ of $\Delta_\alpha$ defined inductively by $z_1 = (x_1 +1) \ \mathrm{mod} \ \alpha (1)$ and
    \begin{equation*}
    z_{i+1}=
        \begin{cases}
     x_{i+1} & \textup{if } z_i \geq x_i \\
    (x_{i+1} +1) \ \mathrm{mod} \ \alpha (i+1) & \textup{if } z_i < x_i .
        \end{cases}
    \end{equation*}
    Clearly, $\varphi_\alpha$ is continuous. The topological dynamical system $(\Delta_\alpha, \varphi_\alpha)$ is called an adding machine or odometer. If it is the case that for every prime number $p$ there are infinitely many $i \in \N$ such that $p$ divides $\alpha (i)$, then $(\Delta_\alpha, \varphi_\alpha)$ is said to be an $\infty$-adic odometer. 
\end{definition}

The following two theorems can be found in \cite{BlockKeesling}.

\begin{theorem}\label{ThmBlockKeesling0}
    For every $\alpha \in (\N \setminus \{ 1 \})^\N$ and every prime number $p$, let $M_\alpha (p)$ be the element of $\N_0 \cup \{ \infty \}$ defined by
    \[ M_\alpha (p) := \sum_{i=1}^\infty \max \big\lbrace n \in \N_0 \, ; \ p^n \textup{ divides } \alpha (i) \big\rbrace . \]
    Then for any $\alpha , \beta \in (\N \setminus \{ 1 \})^\N$, the adding machines $(\Delta_\alpha, \varphi_\alpha)$, $(\Delta_\beta, \varphi_\beta)$ are topologically conjugate if and only if $M_\alpha (p) = M_\beta (p)$ for each prime number $p$. In particular, any two $\infty$-adic odometers are topologically conjugate.
\end{theorem}

\begin{theorem}\label{ThmBlockKeesling}
    Let $\alpha \in (\N \setminus \{ 1 \})^\N$, let $X$ be a compact metric space and let $f \colon X \to X$ be a continuous mapping. Then $(X,f)$ is topologically conjugate to $(\Delta_\alpha, \varphi_\alpha)$ if and only if there is a sequence $(\mathcal{P}_k)_{k=1}^\infty$ of finite clopen partitions of $X$ with $\max \, \lbrace \mathrm{diam} (P) \, ; \, P \in \mathcal{P}_k \rbrace \to 0$ such that, for each $k \in \N$,
    \begin{enumerate}[label=(\arabic*),font=\textup,noitemsep]
        \item $| \mathcal{P}_k | = \prod_{i=1}^k \alpha (i)$;
        \item $\mathcal{P}_{k+1}$ is a refinement of $\mathcal{P}_k$;
        \item members of $\mathcal{P}_k$ are cyclically permuted by $f$, that is, there is an enumeration $P_0, \dotsc ,P_{n-1}$ of members of $\mathcal{P}_k$ such that $f$ maps $P_i$ onto $P_{(i+1) \, \mathrm{mod} \, n}$ for each $i \in \{ 0, \dotsc ,n-1 \}$.
    \end{enumerate}
\end{theorem}

Fix $M>0$ and denote
\begin{align*}
    C(M)&:= \big\lbrace f \colon [-M,M] \to \R \, ; \ f \textup{ is continuous} \big\rbrace ,\\
    C^1(M)&:=\big\lbrace f \in C(M) \, ; \ f \textup{ is continuously differentiable} \big\rbrace .
\end{align*}
Define metrics $d$ and $d_1$ on $C(M)$ and $C^1(M)$, respectively, by
\begin{align*}
    d(f,g) &= \sup \big\lbrace |f(t)-g(t)| \, ; \ t \in [-M,M] \big\rbrace,\\
    d_1(f,g) &= d(f,g)+d(f',g').
\end{align*}
Then $(C(M),d)$ and $(C^1(M),d_1)$ are separable complete metric spaces.

For any $f \in C^1(M)$, we write $f'(M)$ and $f'(-M)$ in place of $f_{-}'(M)$ and $f_{+}'(-M)$, respectively. For all $f \in C^1(M)$ and $x \in [-M,M]$ with $f'(x) \neq 0$, we denote
\[ N(f,x):=x-\frac{f(x)}{f'(x)} . \]
We denote by $\tau (f,x)$ the sequence $(x_n)_{n=0}^{\infty}$ given recursively by $x_0=x$ and by $x_{n+1}=N(f,x_n)$ for each $n \in \N_0$. Of course, $\tau (f,x)$ may not be well-defined. When it is well-defined, we let
\[ \omega (f,x):= \bigcap_{m=0}^{\infty} \overline{ \big\lbrace \tau (f,x)(n) \, ; \, n \geq m \big\rbrace } ,\]
where $\tau (f,x)(n)=x_n$ is the $n$-th term of $\tau (f,x)$ for each $n \in \N_0$.

\begin{definition}
    Let $f \in C^1(M)$ and $x \in [-M,M]$. If there are $n \in \N_0$ and $x_0, \dotsc ,x_n \in [-M,M]$ with $x_0=x$ such that $x_{i+1}=N(f,x_i)$ for each $i \in \{0, \dotsc ,n-1 \}$ and $|N(f,x_n)|>M$, we say that Newton's method diverges for $(f,x)$.

\end{definition}

For every $f \in C^1 (M)$, let
\[ \mathrm{Aff}(f) := \big\lbrace J \subseteq [-M,M] \, ; \ J \textup{ is a nondegenerate closed interval, } f\restriction_J \textup{ is affine and nonconstant} \big\rbrace . \]
For any $J \in \mathrm{Aff}(f)$, we denote by $N(f,J)$ the unique $y \in \R$ satisfying $N(f,x)=y$ for every $x \in J$.

\begin{definition}
    Let $f \in C^1 (M)$ and let $\mathcal{F} \subseteq \mathrm{Aff}(f)$ be a finite disjoint family.
    \begin{itemize}
        \item We say that $f$ is \textit{nice} to $\mathcal{F}$ provided that for every $J \in \mathcal{F}$, either $N(f,J) \notin [-M,M]$, or $N(f,J) \in \Int (I)$ for some $I \in \mathcal{F}$.
        \item We say that $\mathcal{F}$ is \textit{cyclic} with respect to $f$ if there is $n \in \N$ and pairwise distinct intervals $J_0, \dotsc ,J_n \in \mathcal{F}$ such that $\{ J_0, \dotsc ,J_n \} = \mathcal{F}$, $N(f,J_n) \in \Int (J_0)$ and $N(f,J_i) \in \Int (J_{i+1})$ for $i \in \{ 0, \dotsc ,n-1 \}$.
    \end{itemize}
\end{definition}

Note that if $f$ is nice to $\mathcal{F}$, then, for any $J \in \mathcal{F}$, either Newton's method diverges for $(f,x)$ for each $x \in J$, or $\tau (f,x)$ is eventually periodic for each $x \in J$.

\section{The results}

\begin{observation}\label{Obser}
    Let $f \in C^1(M)$ and $J \in \mathrm{Aff}(f)$. Then, given any $\varepsilon > 0$, there exists $\eta > 0$ such that $|N(g,x) - N(f,J)|< \varepsilon$ for every $x \in J$ and every $g \in C^1(M)$ with $d_1(f,g)< \eta$.
\end{observation}

\begin{proposition}\label{PropositionConvergenceOfNewtonsMethod}
    Let $f \in C^1 (M)$, $J \in \mathrm{Aff}(f)$ and assume that $f(z)=0$ for some $z \in \Int (J)$. Then there exists $\delta >0$ such that for every $x \in J$ and every $g \in C^1(M)$ with $d_1(f,g)<\delta$, the sequence $\tau (g,x)$ is well-defined, it converges to a root of $g$ and all of its members belong to $J$.
\end{proposition}

\begin{proof}
    There are $a,b \in \R$ with $a \neq 0$ such that $f(x)=ax+b$ for every $x \in J$. We can assume without loss of generality that $a>0$. Clearly, $z=-b/a$. Fix $\varepsilon > 0$ such that $[z-\varepsilon, z+\varepsilon] \subseteq J$ and $\varepsilon < 1/5$. By Observation \ref{Obser}, there is $\eta > 0$ such that $|N(g,x) - N(f,J)|< \varepsilon$ for all $x \in J$ and $g \in C^1(M)$ with $d_1(f,g)< \eta$. Define $\delta := \min \{ \eta , a \varepsilon \}$ and let $g \in C^1(M)$ with $d_1(f,g)<\delta$ be given. Since $\delta \leq \eta$, we have $|N(g,x)-z|=|N(g,x)-N(f,J)|< \varepsilon$ for every $x \in J$. This shows that $\tau (g,x)$ is well-defined and contained in $J$ for every $x \in J$. Note that $g'(x) > f'(x)- \delta \geq a - a \varepsilon > 4a/5$ for every $x \in J$. In particular, $g$ is strictly increasing on $J$. Moreover,
    \[ g(z-\varepsilon) < f(z-\varepsilon) +\delta \leq f(z-\varepsilon) + a \varepsilon = az-a\varepsilon+b+a\varepsilon =0 \]
    and
    \[ g(z+\varepsilon) > f(z+\varepsilon) -\delta \geq f(z+\varepsilon) - a \varepsilon = az+a\varepsilon+b-a\varepsilon =0. \]
    Hence, $g$ has exactly one root $w$ in $J$ and this root satisfies $z-\varepsilon < w < z+\varepsilon$.

    \begin{claim}\label{claim}
        For every $x \in J$, we have $|w-N(g,x)| \leq |w-x|/2$.
    \end{claim}
    
    \begin{claimproof}
        Fix $x \in J$ and let $\widehat{x}:=N(g,x)$. If $x=w$, then $\widehat{x}=N(g,w)=w-0=w$ and we are done. Assume $x<w$. By the mean value theorem, there is $\xi \in (x,w)$ such that
        \[ g'(\xi)=\frac{g(w)-g(x)}{w-x}=\frac{-g(x)}{w-x}=\frac{(\widehat{x}-x)g'(x)}{w-x}. \]
        Since $a-\delta<g'(\xi)<a+\delta$ and both $g'(x)$ and $(w-x)$ are positive numbers, we have
        \[ \frac{a-\delta}{g'(x)}(w-x)<\widehat{x}-x<\frac{a+\delta}{g'(x)}(w-x). \]
        Subtracting $(w-x)$ we obtain
        \[ \frac{a-\delta -g'(x)}{g'(x)}(w-x)<\widehat{x}-w<\frac{a+\delta -g'(x)}{g'(x)}(w-x). \]
        Since $0 < a-\delta<g'(x) < a+\delta$, this leads to
        \[ \frac{a-\delta -(a+\delta)}{a-\delta}(w-x)<\widehat{x}-w<\frac{a+\delta -(a-\delta)}{a-\delta}(w-x), \]
        hence
        \[ \frac{-2\delta}{a-\delta}(w-x)<\widehat{x}-w<\frac{2\delta}{a-\delta}(w-x). \]
        As $\delta \leq a \varepsilon < a/5$, it follows that
        \[ \frac{-2 \cdot (a/5)}{a-(a/5)}(w-x) < \widehat{x}-w < \frac{2 \cdot (a/5)}{a-(a/5)}(w-x), \]
        thus $(x-w)/2<\widehat{x}-w<(w-x)/2$, which proves that $|w-\widehat{x}| \leq |w-x|/2$. A similar approach works if $x>w$. \claimend
    \end{claimproof}
    Finally, given any $x \in J$, it follows from Claim \ref{claim} that $|w-\tau (g,x)(n)|\leq 2^{-n} |w-x|$ for each $n \in \N$, hence $\tau (g,x)$ converges to $w$.
\end{proof}

\begin{proposition}\label{PropositionApproximationByOftenAffineFunctions1}
    Let $f \in C^1(M)$, $\varepsilon > 0$, $\delta > 0$ and $t \in (-M,M)$. There exists $g \in C^1(M)$ together with a finite disjoint family $\mathcal{F} \subseteq \mathrm{Aff}(g)$ such that:
    \begin{enumerate}[label=(\arabic*),font=\textup,noitemsep]
        \item $d_1(f,g)< \varepsilon$;
        \item $\lambda \big( \bigcup \mathcal{F} \big) > 2M-\delta$;
        \item $t \in \Int \big( \bigcup \mathcal{F} \big)$;
        \item $g$ is nice to $\mathcal{F}$.
    \end{enumerate}
\end{proposition}

\begin{proof}
    Since $f$ and $f'$ are continuous on $[-M,M]$, there is $r>0$ such that $|f(x)-f(y)|< \frac{1}{12} \varepsilon$ and $|f'(x)-f'(y)|< \frac{1}{12} \varepsilon$ for all $x,y \in [-M,M]$ with $|x-y|<2r$. Let $z_0, z_1, \dotsc ,z_n \in \R \setminus \{ t \}$ satisfy $-M = z_0 < z_1 < \dotsc < z_n = M$ and $|z_i - z_{i-1}| < r$ for every $i \in \{ 1, \dotsc ,n \}$. Let $w_0, \dotsc ,w_n$ be pairwise distinct real numbers such that:
    \begin{enumerate}[label=(\roman*),font=\textup]
        \item $|w_i - f(z_i)|< \varepsilon /12 \,$ for each $i \in \{ 0, \dotsc ,n \}$;
        \item $\displaystyle \bigg| \frac{f(z_i)-f(z_{i-1})}{z_i-z_{i-1}} - \frac{w_i-w_{i-1}}{z_i-z_{i-1}} \bigg| < \varepsilon /12 \,$ for each $i \in \{ 1, \dotsc ,n \}$;
        \item for all $i \in \{ 1, \dotsc ,n \}$, the root $x_i$ of $\displaystyle s \mapsto \frac{w_i-w_{i-1}}{z_i-z_{i-1}} (s-z_i) + w_i$ does not belong to $\{ z_0, \dotsc ,z_n \}$.
    \end{enumerate}
    Note that $|w_i-w_{i-1}| \leq |w_i-f(z_i)| + |f(z_i)-f(z_{i-1})| + |f(z_{i-1}) - w_{i-1}| < 3 \varepsilon /12 = \varepsilon /4$ for each $i \in \{ 1, \dotsc ,n \}$. Define a function $h \colon [-M,M] \to \R$ by
    \[ h (s) = \frac{w_i-w_{i-1}}{z_i-z_{i-1}} (s-z_i) + w_i \, , \ \ i \in \{ 1, \dotsc ,n \} , \, s \in [z_{i-1},z_i]. \]
    Then $h \in C(M)$ and, for all $i \in \{ 1, \dotsc ,n \}$ and $s \in [z_{i-1},z_i]$, we have
    \[ |w_i - h(s)| = |h(z_i) - h(s)| \leq |h(z_i) - h(z_{i-1})| = |w_i-w_{i-1}| < \varepsilon /4 , \]
    hence
    \[ |f(s)-h(s)| \leq |f(s)-f(z_i)| + |f(z_i)-w_i| + |w_i-h(s)| < (2 \varepsilon /12) + (\varepsilon /4) < \varepsilon /2 . \]
    This proves that $d(f,h)< \varepsilon /2$. For every $i \in \{ 1, \dotsc ,n \}$, by Lagrange's mean value theorem, there is $\xi_i \in (z_{i-1},z_i)$ such that $f(z_i)-f(z_{i-1}) = (z_i-z_{i-1})  f'(\xi_i)$. Hence, for every $s \in (z_{i-1},z_i)$,
    \[ |f'(s)-h'(s)| \leq |f'(s)-f'(\xi_i)| + \bigg| \frac{f(z_i)-f(z_{i-1})}{z_i-z_{i-1}} - \frac{w_i-w_{i-1}}{z_i-z_{i-1}} \bigg| < 2 \varepsilon /12 = \varepsilon /6 .\]
    Let $V_1, \dotsc ,V_{n-1}$ be open subintervals of $(-M,M)$ such that
    \begin{enumerate}[label=(\alph*),noitemsep]
        \item $z_i \in V_i$ for every $i \in \{ 1, \dotsc ,n-1 \}$;
        \item $\overline{V_i} \cap \overline{V_j} = \emptyset$ for all $i,j \in \{ 1, \dotsc ,n-1 \}$ with $i \neq j$;
        \item $\overline{V_i} \cap \big( \{ -M,M,t \} \cup \{ x_1, \dotsc ,x_n \} \big) = \emptyset$ for every $i \in \{ 1, \dotsc ,n-1 \}$;
        \item $\lambda (V_1 \cup \dots \cup V_{n-1}) < \delta$.
    \end{enumerate}
    Define
    \[ \mathcal{F}:= \big\lbrace K \subseteq [-M,M] \, ; \ K \textup{ is a connected component of } [-M,M] \setminus (V_1 \cup \dots \cup V_{n-1}) \big\rbrace . \]
    Clearly, $\mathcal{F}$ is a finite family of pairwise disjoint nondegenerated closed subintervals of $[-M,M]$. Moreover, $t \in \Int \big( \bigcup \mathcal{F} \big)$ and
    \[ \textstyle \lambda \big( \bigcup \mathcal{F} \big) = \lambda \big( [-M,M] \setminus (V_1 \cup \dots \cup V_{n-1}) \big) = 2M - \lambda (V_1 \cup \dots \cup V_{n-1}) > 2M - \delta . \]
    It is not difficult to construct (see e.g. \cite{Misiurewicz}) a function $g \in C^1(M)$ such that
    \begin{itemize}[noitemsep]
        \item $d(g,h) < \varepsilon /12 \,$;
        \item $g(s)=h(s)$ for every $s \in \bigcup \mathcal{F}$;
        \item $\min \big\lbrace h_{-}'(z_i), h_{+}'(z_i) \big\rbrace \leq g'(s) \leq \max \big\lbrace h_{-}'(z_i), h_{+}'(z_i) \big\rbrace$ for all $s \in V_i$ and $i \in \{ 1, \dotsc ,n-1 \}$.
    \end{itemize}
    Then $\mathcal{F} \subseteq \mathrm{Aff}(g)$ and it is easy to see that $g$ is nice to $\mathcal{F}$. For all $i \in \{ 1, \dotsc ,n-1 \}$ and $s \in V_i \setminus \{ z_i \}$, we have $h'(s) \in \big\lbrace h_{-}'(z_i), h_{+}'(z_i) \big\rbrace$, thus
    \begin{align*}
        |h'(s)-g'(s)| &\leq \big| h_{-}'(z_i) - h_{+}'(z_i) \big| = \bigg| \frac{w_i-w_{i-1}}{z_i-z_{i-1}} - \frac{w_{i+1}-w_i}{z_{i+1}-z_i} \bigg| \\
        &\leq \bigg| \frac{w_i-w_{i-1}}{z_i-z_{i-1}} - f'(\xi_i) \bigg| + \big| f'(\xi_i) - f'(\xi_{i+1}) \big| + \bigg| f'(\xi_{i+1}) - \frac{w_{i+1}-w_i}{z_{i+1}-z_i} \bigg| < \frac{3 \varepsilon}{12} = \frac{\varepsilon}{4} .
    \end{align*}
    Consequently, $|f'(s)-g'(s)| \leq |f'(s)-h'(s)| + |h'(s)-g'(s)| < (\varepsilon /6) + (\varepsilon /4) = 5 \varepsilon /12$ for every $s \in [-M,M] \setminus \{ z_0, \dotsc ,z_n \}$, hence
    \[ d_1(f,g) = d(f,g) + d(f',g') \leq d(f,h) + d(h,g) + d(f',g') < (\varepsilon /2) + (\varepsilon /12) + (5 \varepsilon /12) = \varepsilon . \]
\end{proof}

\begin{lemma}\label{Lemma}
    Let $f \in C^1(M)$, $J \in \mathrm{Aff}(f)$ and $\varepsilon > 0$. Assume that $N(f,J) \notin J$. Then there exists $r>0$ such that for every 
    $\delta > 0$, every $x \in \Int (J)$ and every $y \in \R$ with $|y-N(f,J)|<r$, there is $g \in C^1(M)$ satisfying:
    \begin{enumerate}[label=(\roman*),font=\textup,noitemsep]
        \item $d_1(f,g)<\varepsilon$;
        \item $f(t)=g(t)$ for every $t \in [-M,M] \setminus (x-\delta,x+\delta)$;
        \item $g$ is affine and nonconstant on some neighbourhood of $x$;
        \item $N(g,x)=y$.
    \end{enumerate}
\end{lemma}

\begin{proof}
    As $J \in \mathrm{Aff}(f)$, there are $a,b \in \R$ with $a \neq 0$ such that $f(t)=at+b$ for every $t \in J$. Denote $z := N(f,J)$. Clearly, $z = -b/a$. Since $z \notin J$ and $J$ is closed, we have $\mathrm{dist}(z,J)>0$. Take $r \in \R$ with $0 < r < \mathrm{dist}(z,J)$ such that
    \[ \frac{|a|r}{\mathrm{dist}(z,J)-r} < \frac{\varepsilon}{2}. \]
    Let $\delta > 0$, $x \in \Int (J)$ and $y \in (z-r,z+r)$ be given. Fix $\gamma \in (0,1]$ such that $[x-2 \gamma, x+2 \gamma] \subseteq \Int(J)$ and $2 \gamma < \delta$. Let $g_0 \colon [-M,M] \to \R$ be the unique continuous function satisfying
        \begin{itemize}
        \item $g_0(t)=f(t)$ for every $t \in [-M,M] \setminus (x-2\gamma,x+2\gamma)$;
        \item $\displaystyle g_0(t)=\frac{f(x)}{x-y}(t-x)+f(x)$ for every $t \in [x-\gamma,x+\gamma]$;
        \item $g_0$ is affine on both $[x-2\gamma,x-\gamma]$ and $[x+\gamma,x+2\gamma]$.
    \end{itemize}
    Then $N(g_0,x)=x- \big( g_0(x)/g_0'(x) \big) = x- \big( f(x)/g_0'(x) \big) = x-(x-y)=y$. Also, since
    \[ \bigg| a - \frac{f(x)}{x-y} \bigg| = \bigg| \frac{ax-ay-ax-b}{x-y} \bigg| = \bigg| \frac{az-ay}{x-y} \bigg| < \frac{|a|r}{\mathrm{dist}(y,J)} \leq \frac{|a|r}{\mathrm{dist}(z,J)-r} < \frac{\varepsilon}{2}, \]
    it follows that $|f'(t)-g_0'(t)|< \varepsilon /2$ for every $t \in [-M,M] \setminus \{ x - 2\gamma, \, x - \gamma, \, x + \gamma, \, x +2 \gamma \}$. Moreover, it is easy to see that
    \[ d(f,g_0) = \big| g_0 (x + \gamma) - f(x + \gamma) \big| = \bigg| \frac{\gamma f(x)}{x-y}+f(x)- \big( f(x) + a \gamma \big) \bigg| = \gamma \bigg| \frac{f(x)}{x-y} -a \bigg| < \frac{\varepsilon}{2} . \]
    Redefining $g_0$ on sufficiently small neighborhoods of $x-2\gamma$, $x-\gamma$, $x+\gamma$ and $x+2\gamma$, we obtain the desired function $g \in C^1(M)$.
\end{proof}

\begin{proposition}\label{PropositionChangingThePeriodOfACycle}
    Let $f \in C^1(M)$ and let $\mathcal{F} \subseteq \mathrm{Aff}(f)$ be a finite disjoint family such that $f$ is nice to $\mathcal{F}$. Let $\mathcal{J} \subseteq \mathcal{F}$ be cyclic with respect to $f$, let $\varepsilon>0$, $\Delta>0$, $m \in \N$ and $t \in \Int \big( \bigcup \mathcal{J} \big)$. There is $g \in C^1(M)$ and a finite disjoint family $\mathcal{S} \subseteq \mathrm{Aff}(g)$ such that:
    \begin{enumerate}[label=(\arabic*),font=\textup,noitemsep]
        \item $\bigcup \mathcal{S} \subseteq \bigcup \mathcal{J}$;
        \item $f(x)=g(x)$ for every $x \in [-M,M] \setminus \bigcup \mathcal{J}$;
        \item $\mathcal{F} \setminus \mathcal{J} \subseteq \mathrm{Aff}(g)$ and $g$ is nice to $(\mathcal{F} \setminus \mathcal{J}) \cup \mathcal{S}$;
        \item $d_1(f,g)< \varepsilon$;
        \item $\lambda \big( \big( \bigcup \mathcal{J} \big) \setminus \big( \bigcup \mathcal{S} \big) \big) < \Delta$;
        \item $t \in \Int \big( \bigcup \mathcal{S} \big)$;
        \item for each $x \in \bigcup \mathcal{S}$, $\tau (g,x)$ is eventually periodic with period greater than $1$ and divisible by $m$.
    \end{enumerate}
\end{proposition}

\begin{proof}
    There is $n \in \N$ and pairwise distinct intervals $J_0, \dotsc ,J_n \in \mathcal{J}$ such that $\{ J_0, \dotsc ,J_n \} = \mathcal{J}$, $z_0 := N(f,J_n) \in \Int (J_0)$ and $z_i := N(f,J_{i-1}) \in \Int (J_i)$ for $i \in \{ 1, \dotsc ,n \}$. Let
    \[ F := \big( \{ t \} \cup \{ N(f,I) \, ; \, I \in \mathcal{F} \} \big) \setminus \{ z_0, \dotsc ,z_n \} . \]
    For each $i \in \{ 0, \dotsc ,n \}$, take $r_i > 0$ from Lemma \ref{Lemma} for $f$, $J_i$ and $\varepsilon$. For every $i \in \{ 0, \dotsc ,n \}$, fix $m$ distinct points $x_1^i, \dotsc ,x_m^i \in \Int (J_i) \setminus F$ such that $z_i \in \{ x_1^i, \dotsc ,x_m^i \}$ and $|z_i - x_j^i| < r_i$ for each $j \in \{ 1, \dotsc ,m \}$. Fix $\delta > 0$ such that, for every $i \in \{ 0, \dotsc ,n \}$,
    \begin{itemize}
        \item $[x_1^i-\delta, x_1^i+\delta], \dotsc ,[x_m^i-\delta, x_m^i+\delta]$ are pairwise disjoint intervals contained in $\Int (J_i)$;
        \item $[x_j^i-\delta, x_j^i+\delta] \cap F = \emptyset$ for every $j \in \{ 1, \dotsc ,m \}$;
        \item $2m\delta(n+1)< \Delta$.
    \end{itemize}
    For all $i \in \{ 0, \dotsc ,n \}$ and $j \in \{ 1, \dotsc ,m \}$, there is $g_j^i \in C^1(M)$ such that:
    \begin{enumerate}[label=(\roman*),font=\textup,noitemsep]
        \item $d_1(f,g_j^i)<\varepsilon$;
        \item $f(s)=g_j^i(s)$ for every $s \in [-M,M] \setminus (x_j^i-\delta,x_j^i+\delta)$;
        \item $g_j^i$ is affine and nonconstant on some neighbourhood of $x_j^i$;
    \end{enumerate}
    and
    \begin{equation*}
        N(g_j^i,x_j^i)=
        \begin{cases}
            x_j^{i+1} & \textup{if } i < n \\
            x_{j+1}^0 & \textup{if } i=n \textup{ and } j<m \\
            x_1^0 & \textup{if } i=n \textup{ and } j=m.
        \end{cases}
    \end{equation*}
    For all $i \in \{ 0, \dotsc ,n \}$ and $j \in \{ 1, \dotsc ,m \}$, let $V_j^i := (x_j^i-\delta,x_j^i+\delta)$ and fix a closed interval $K_j^i \subseteq V_j^i$ with $x_j^i \in \Int (K_j^i)$ such that $g_j^i$ is affine on $K_j^i$. Denote $G:= \bigcup \big\lbrace V_j^i \, ; \ i \in \{ 0, \dotsc ,n \}, j \in \{ 1, \dotsc ,m \} \big\rbrace$ and $\mathcal{K}:=\big\lbrace K_j^i \, ; \ i \in \{ 0, \dotsc ,n \}, j \in \{ 1, \dotsc ,m \} \big\rbrace$. Finally, let
    \[ \mathcal{S} := \mathcal{K} \cup \bigcup_{i=0}^n \big\lbrace L \subseteq J_i \setminus G \, ; \ L \textup{ is a connected component of } J_i \setminus G \big\rbrace \]
    and define a function $g \colon [-M,M] \to \R$ by $g(s)=f(s)$ for $s \in [-M,M] \setminus G$ and by $g(s)=g_j^i(s)$ for $s \in V_j^i$, $i \in \{ 0, \dotsc ,n \}$ and $j \in \{ 1, \dotsc ,m \}$. Clearly, $(\mathcal{F} \setminus \mathcal{J}) \cup \mathcal{S} \subseteq \mathrm{Aff}(g)$. Since $\overline{G} \cap F = \emptyset$ and $z_i \in \{ x_1^i, \dotsc ,x_m^i \} \subseteq \Int (K_1^i) \cup \dots \cup \Int (K_m^i) \subseteq \bigcup \{ \Int (S) \, ; \, S \in \mathcal{S} \}$ for each $i \in \{ 0, \dotsc ,n \}$, it easily follows that $g$ is nice to $(\mathcal{F} \setminus \mathcal{J}) \cup \mathcal{S}$ and that $t \in \Int \big( \bigcup \mathcal{S} \big)$. Obviously, conditions (1), (2) and (4) are satisfied as well. For each $x \in \bigcup \mathcal{S}$, it follows from the construction that $\tau (g,x)$ is eventually periodic with period $m(n+1)$. Moreover, since $G \subseteq \bigcup \mathcal{J}$ and $(\bigcup \mathcal{J}) \setminus G \subseteq \bigcup \mathcal{S}$, we obtain
    \[ \textstyle \lambda \Big( \big( \bigcup \mathcal{J} \big) \setminus \big( \bigcup \mathcal{S} \big) \Big) \leq \lambda \Big( \big( \bigcup \mathcal{J} \big) \setminus \big( \big( \bigcup \mathcal{J} \big) \setminus G \big) \Big) = \lambda (G) = \displaystyle \sum_{i=0}^n \sum_{j=1}^m \lambda (V_j^i) = 2m\delta(n+1)< \Delta . \]
\end{proof}

By Propositions \ref{PropositionApproximationByOftenAffineFunctions1} and \ref{PropositionChangingThePeriodOfACycle}, we have the following:

\begin{corollary}\label{CorolA}
    Let $f \in C^1(M)$, $\varepsilon > 0$, $\delta > 0$, $m \in \N$ and $t \in (-M,M)$. There exists $g \in C^1(M)$ together with a finite disjoint family $\mathcal{A} \subseteq \mathrm{Aff}(g)$ such that:
    \begin{enumerate}[label=(\arabic*),font=\textup,noitemsep]
        \item $g$ is nice to $\mathcal{A}$;
        \item $d_1(f,g)< \varepsilon$;
        \item $\lambda \big( \bigcup \mathcal{A} \big) > 2M-\delta$;
        \item $t \in \Int \big( \bigcup \mathcal{A} \big)$;
        \item for every $J \in \mathcal{A}$, either Newton's method diverges for $(g,x)$ for each $x \in J$, or $\tau (g,x)$ is eventually constant for each $x \in J$, or $\tau (g,x)$ is eventually periodic with period greater than $1$ and divisible by $m$ for each $x \in J$.
    \end{enumerate}
\end{corollary}

\begin{proposition}\label{Prop}
    Let $\delta >0$. There is a dense $G_{\delta}$ set $\mathcal{G} \subseteq C^1(M)$ such that for every $h \in \mathcal{G}$ there is a closed set $S \subseteq [-M,M]$ satisfying the following two conditions:
    \begin{enumerate}[label=(\arabic*),font=\textup,noitemsep]
        \item $\lambda (S) > 2M-\delta$;
        \item for all $x \in S$, either Newton's method diverges for $(h,x)$, or 
$\tau (h,x)$ converges to a root of $h$, or $\omega (h,x)$ is an $\infty$-adic odometer.
    \end{enumerate}
\end{proposition}

\begin{proof}
    Let $\{ f_i \, ; \, i \in \N \}$ be a dense subset of $C^1(M)$. For all $i,j \in \N$, there is (by Corollary \ref{CorolA}) a function $g_{i,j} \in C^1(M)$ together with a finite disjoint family $\mathcal{A}_{i,j} \subseteq \mathrm{Aff}(g_{i,j})$ such that:
    \begin{enumerate}[label=(\arabic*),font=\textup,noitemsep]
        \item $g_{i,j}$ is nice to $\mathcal{A}_{i,j}$;
        \item $d_1(f_i,g_{i,j})< 2^{-i-j} \delta$;
        \item $\lambda \big( \bigcup \mathcal{A}_{i,j} \big) > 2M-2^{-i-j} \delta$;
        \item for every $J \in \mathcal{A}_{i,j}$, either Newton's method diverges for $(g_{i,j},x)$ for each $x \in J$, or $\tau (g_{i,j},x)$ is eventually constant for each $x \in J$, or $\tau (g_{i,j},x)$ is eventually periodic with period greater than $1$ and divisible by $j$ for each $x \in J$.
    \end{enumerate}
    For all $i,j \in \N$, let
    \begin{align*}
        \mathcal{F}_{i,j}^1 &:= \big\lbrace J \in \mathcal{A}_{i,j} \, ; \ |N(g_{i,j},J)| > M \big\rbrace ,\\
        \mathcal{F}_{i,j}^2 &:= \big\lbrace J \in \mathcal{A}_{i,j} \, ; \ N(g_{i,j},J) \in \Int (J) \big\rbrace ,\\
        \mathcal{F}_{i,j}^3 &:= \big\lbrace J \in \mathcal{A}_{i,j} \, ; \ N(g_{i,j},J) \in \Int (I) \textup{ for some } I \in \mathcal{A}_{i,j} \setminus \{ J \} \big\rbrace .
    \end{align*}
    Note that the families $\mathcal{F}_{i,j}^1, \mathcal{F}_{i,j}^2, \mathcal{F}_{i,j}^3$ are pairwise disjoint and that (as $g_{i,j}$ is nice to $\mathcal{A}_{i,j}$) their union is equal to $\mathcal{A}_{i,j}$. For all $i,j \in \N$ and $J \in \mathcal{F}_{i,j}^3$, denote by $I_{i,j,J}$ the member of $\mathcal{A}_{i,j} \setminus \{ J \}$ satisfying $N(g_{i,j},J) \in \Int (I_{i,j,J})$.
    
    For all $i,j \in \N$ and $J \in \mathcal{F}_{i,j}^1$, it follows from Observation \ref{Obser} that there exists $\eta_{i,j,J}>0$ such that $|N(h,y)|>M$ for every $y \in J$ and every $h \in C^1(M)$ with $d_1(h,g_{i,j})< \eta_{i,j,J}$. By Proposition \ref{PropositionConvergenceOfNewtonsMethod}, for all $i,j \in \N$ and $J \in \mathcal{F}_{i,j}^2$, there is $\eta_{i,j,J}>0$ such that for every $y \in J$ and every $h \in C^1(M)$ with $d_1(h,g_{i,j})<\eta_{i,j,J}$, the sequence $\tau (h,y)$ converges to a root of $h$. Finally, for all $i,j \in \N$ and $J \in \mathcal{F}_{i,j}^3$, it follows from Observation \ref{Obser} that there are $\eta_{i,j,J}>0$ and a nondegenerate closed interval of length less than $1/j$ contained in $I_{i,j,J}$ such that $N(h,y)$ belongs to that interval for all $y \in J$ and $h \in C^1(M)$ with $d_1(h,g_{i,j})< \eta_{i,j,J}$. 

    Let $\eta_{i,j} := \min \big( \{ 2^{-i-j} \delta \} \cup \{ \eta_{i,j,J} \, ; \, J \in \mathcal{A}_{i,j} \} \big)$ and $\mathcal{B}_{i,j} := \big\lbrace h \in C^1(M) \, ; \ d_1(h,g_{i,j}) < \eta_{i,j} \big\rbrace$ for all $i,j \in \N$. Finally, let
    \[ \mathcal{G} := \bigcap_{j \in \N} \bigcup_{i \in \N} \mathcal{B}_{i,j} . \]
    Then $\mathcal{G}$ is $G_\delta$ in $C^1(M)$ and, by the Baire category theorem, it is dense. Fix $h \in \mathcal{G}$. For each $j \in \N$, there is $i(j) \in \N$ such that $h \in \mathcal{B}_{i(j),j}$. Let $S := \bigcap \big\lbrace \bigcup \mathcal{A}_{i(j),j} \, ; \ j \in \N \big\rbrace$. Then $S$ is closed and
    \[ \lambda \big( [-M,M] \setminus S \big) \leq \sum_{j \in \N} \lambda \big( [-M,M] \setminus \textstyle \bigcup \mathcal{A}_{i(j),j} \big) < \displaystyle \sum_{j \in \N} 2^{-i(j)-j} \delta \leq \delta , \]
    hence $\lambda (S) > 2M-\delta$. Fix $x \in S$. If $\tau (g_{i(j),j},x)$ is not well-defined for some $j \in \N$, then Newton's method for $(g_{i(j),j},x)$ diverges (as $x$ belongs to a member of $\mathcal{A}_{i(j),j}$ and $g_{i(j),j}$ is nice to $\mathcal{A}_{i(j),j}$) and so does Newton's method for $(h,x)$ due to the fact that $d_1 (h,g_{i(j),j}) < \eta_{i(j),j}$. Thus, for every $j \in \N$, we can assume that $\tau (g_{i(j),j},x)$ is well-defined and either eventually constant, or eventually periodic with period greater than $1$ and divisible by $j$. If $\tau (g_{i(j),j},x)$ is eventually constant for some $j \in \N$, then, since $d_1 (h,g_{i(j),j}) < \eta_{i(j),j}$ and $x$ belongs to a member of $\mathcal{A}_{i(j),j}$, it follows that $\tau (h,x)$ converges to a root of $h$. Finally, assume that for every $j \in \N$, $\tau (g_{i(j),j},x)$ is eventually periodic with period greater than $1$ and divisible by $j$. Then, for each $j \in \N$, since $d_1 (h,g_{i(j),j}) < \eta_{i(j),j}$, there is a natural number $n(j)$ divisible by $j$ and pairwise disjoint nondegenerate closed intervals $J_1^j, \dotsc ,J_{n(j)}^j$ contained in $[-M,M]$ such that:
    \begin{itemize}[noitemsep]
        \item $\lambda (J_k^j) < 1/j$ for each $k \in \{ 1, \dotsc ,n(j) \}$;
        \item $N(h,y) \in J_{k+1}^j$ for all $y \in J_k^j$ and $k \in \{ 1, \dotsc ,n(j)-1 \}$;
        \item $N(h,y) \in J_1^j$ for all $y \in J_{n(j)}^j$;
        \item $\tau (h,x)(k) \in J_1^j \cup \dots \cup J_{n(j)}^j$ for some $k \in \N_0$.
    \end{itemize}
    For every $j \in \N$, let $\mathcal{P}_j := \{ J_k^j \cap \omega (h,x) \, ; \ k = 1, \dotsc ,n(j) \} $. Then $\mathcal{P}_j$ is a partition of $\omega (h,x)$ consisting of relatively clopen sets cyclically permuted by $N(h, \cdot)$. Let $(j_k)_{k=1}^\infty$ be an increasing sequence of natural numbers greater than $1$ such that, for each $k \in \N$,
    \begin{itemize}[noitemsep]
        \item $j_{k+1}$ is divisible by $k \cdot n(j_k)$;
        \item $1/j_{k+1} < \min \big\lbrace \mathrm{dist}(I,J) \, ; \ I,J \textup{ are distinct members of } \mathcal{P}_{j_k} \big\rbrace$.
    \end{itemize}
    Let $\alpha (1) := n(j_1)$ and $\alpha (k+1) := n(j_{k+1}) / n(j_k)$ for each $k \in \N$. Then, for every $k \in \N$, $\alpha (k)$ is a natural number greater than $1$ and it is divisible by $k$. Letting $\alpha := ( \alpha (1), \alpha (2), \dotsc )$, it follows that $(\Delta_\alpha, \varphi_\alpha)$ is an $\infty$-adic odometer. For each $k \in \N$, it is easy to see that $\mathcal{P}_{j_{k+1}}$ is a refinement of $\mathcal{P}_{j_k}$ and that $|\mathcal{P}_{j_k}|=n(j_k)= \alpha (1) \cdot \dotsc \cdot \alpha (k)$. Thus, $\big( \omega (h,x), N(h,\cdot) \big)$ is topologically conjugate to an $\infty$-adic odometer by Theorem \ref{ThmBlockKeesling}.
\end{proof}

\begin{corollary}\label{CorolB}
    There exists a dense $G_\delta$ set $\mathcal{G} \subseteq C^1(M)$ such that for any $h \in \mathcal{G}$, there is an $F_\sigma$ set $S \subseteq [-M,M]$ satisfying the following two conditions:
    \begin{enumerate}[label=(\arabic*),font=\textup,noitemsep]
        \item $\lambda (S) = 2M$;
        \item for all $x \in S$, either Newton's method diverges for $(h,x)$, or $\tau (h,x)$ converges to a root of $h$, or $\omega (h,x)$ is an $\infty$-adic odometer.
    \end{enumerate}
\end{corollary}

\begin{theorem}\label{ThmA}
    There is a dense $G_\delta$ set $\mathcal{H} \subseteq C^1(M) \times [-M,M]$ such that for any $(h,x) \in \mathcal{H}$, either Newton's method diverges for $(h,x)$, or $\tau (h,x)$ converges to a root of $h$, or $\omega (h,x)$ is an $\infty$-adic odometer.
\end{theorem}

\begin{proof}
    Let $\{ (f_i,t_i) \, ; \ i \in \N \}$ be a dense subset of $C^1(M) \times (-M,M)$. For all $i,j \in \N$, by Corollary \ref{CorolA}, there is $g_{i,j} \in C^1(M)$ together with a finite disjoint family $\mathcal{A}_{i,j} \subseteq \mathrm{Aff}(g_{i,j})$ such that:
    \begin{enumerate}[label=(\arabic*),font=\textup,noitemsep]
        \item $g_{i,j}$ is nice to $\mathcal{A}_{i,j}$;
        \item $d_1(f_i,g_{i,j})< 2^{-i-j}$;
        \item $t_i \in \Int \big( \bigcup \mathcal{A}_{i,j} \big)$;
        \item for every $J \in \mathcal{A}_{i,j}$, either Newton's method diverges for $(g_{i,j},x)$ for each $x \in J$, or $\tau (g_{i,j},x)$ is eventually constant for each $x \in J$, or $\tau (g_{i,j},x)$ is eventually periodic with period greater than $1$ and divisible by $j$ for each $x \in J$.
    \end{enumerate}
    For $i,j \in \N$, define $\eta_{i,j}$ and $\mathcal{B}_{i,j}$ as in the proof of Proposition \ref{Prop} and denote $V_{i,j} := \Int \big( \bigcup \mathcal{A}_{i,j} \big)$. Let
    \[ \mathcal{H} := \bigcap_{j \in \N} \bigcup_{i \in \N} (\mathcal{B}_{i,j} \times V_{i,j}). \]
    Then $\mathcal{H}$ is $G_\delta$ in $C^1(M) \times [-M,M]$ and, thanks to (3) and the Baire category theorem, it is dense. For each $j \in \N$, there is $i(j) \in \N$ such that $(h,x) \in \mathcal{B}_{i(j),j} \times V_{i(j),j}$. Repeating the arguments used in the proof of Proposition \ref{Prop}, we conclude that either Newton's method diverges for $(h,x)$, or $\tau (h,x)$ converges to a root of $h$, or $\omega (h,x)$ is an $\infty$-adic odometer.
\end{proof}

\begin{theorem}\label{ThmB}
    There exists a dense $G_\delta$ set $\mathcal{R} \subseteq C^1(M)$ such that for any $f \in \mathcal{R}$, there is a Borel set $B \subseteq [-M,M]$ satisfying the following three conditions:
    \begin{enumerate}[label=(\arabic*),font=\textup,noitemsep]
        \item $B$ is comeagre in $[-M,M]$;
        \item $\lambda (B) = 2M$;
        \item for all $x \in B$, either Newton's method diverges for $(f,x)$, or $\tau (f,x)$ converges to a root of $f$, or $\omega (f,x)$ is an $\infty$-adic odometer.
    \end{enumerate}
\end{theorem}

\begin{proof}
    Let $\mathcal{G}$ and $\mathcal{H}$ be the sets obtained in Corollary \ref{CorolB} and Theorem \ref{ThmA}, respectively. By the Kuratowski-Ulam theorem \cite[Theorem 8.41]{Kechris}, there is a dense $G_\delta$ set $\mathcal{H}_0 \subseteq C^1(M)$ such that for every $f \in \mathcal{H}_0$, there is a dense $G_\delta$ set $S_f \subseteq [-M,M]$ satisfying $\{ f \} \times S_f \subseteq \mathcal{H}$. Let $\mathcal{R} := \mathcal{G} \cap \mathcal{H}_0$. Then $\mathcal{R}$ is $G_\delta$ in $C^1(M)$ and, by the Baire category theorem, it is dense. Let $f \in \mathcal{R}$ be given. Since $f \in \mathcal{G}$, there is an $F_\sigma$ set $S \subseteq [-M,M]$ with $\lambda (S) = 2M$ such that for every $x \in S$, either Newton's method diverges for $(f,x)$, or $\tau (f,x)$ converges to a root of $f$, or $\omega (f,x)$ is an $\infty$-adic odometer. Let $B := S \cup S_f$.
\end{proof}

\bibliographystyle{alpha}
\bibliography{references}